\documentclass[12pt]{amsart}

\usepackage{amsmath}
\usepackage{amsthm}
\usepackage{amssymb}
\usepackage{amscd}

\newcommand{\gothic}{\mathfrak}

\newcommand{\m}{{\gothic{m}}}

\newcommand{\End}{\operatorname{End{}}}
\newcommand{\Hom}{\operatorname{Hom{}}}
\newcommand{\Ext}{\operatorname{Ext{}}}
\newcommand{\rank}{\operatorname{rank{}}}

\renewcommand{\hat}{\widehat}

\renewcommand{\phi}{\varphi}

\renewcommand{\to}{{\longrightarrow}}

\newcommand{\e}{\operatorname{e}}
\newcommand{\D}{\operatorname{D}}
\newcommand{\Rh}{R^{\text h}}
\newcommand{\Ker}{\operatorname{Ker}}
\newcommand{\Cl}{\operatorname{Cl}}
\newcommand{\Coker}{\operatorname{Coker}}
\newcommand{\Ass}{\operatorname{Ass}}
\newcommand{\Spec}{\operatorname{Spec}}

\newcommand{\N}{\mathbb{N}}

\newtheorem{thm}{Theorem}
\newtheorem{cor}[thm]{Corollary}
\newtheorem{prop}[thm]{Proposition}
\newtheorem{lemma}[thm]{Lemma}

\theoremstyle{definition}

\newtheoremstyle{newremark}{}{}{\textrm}{}{\textbf}{}{}{}

\newtheorem{eg}[thm]{Example}

\begin{document}

\title{Local Rings of Bounded Cohen--Macaulay Type}

\author{Graham J. Leuschke}

\address{Department of Mathematics \\
        University of Kansas \\
        Lawrence, KS
        66045}

\email{gleuschke@math.ukans.edu}

\urladdr{http://www.leuschke.org}

\author{Roger Wiegand}

\address{Department of Mathematics and Statistics\\
     University of Nebraska--Lincoln\\
    Lincoln, NE
    68588-0323 }

\email{rwiegand@math.unl.edu}
\urladdr{http://www.math.unl.edu/\textasciitilde rwiegand}

\date{April 15, 2003 (filename: BCMTrevised1.tex)}

\thanks{Leuschke's research was supported by an NSF Postdoctoral Fellowship, and Wiegand's was supported by grants from the NSA and the NSF}

\numberwithin{thm}{section}
\renewcommand{\thesubsection}{thm}
\numberwithin{equation}{thm}

\begin{abstract}  Let $(R,\m,k)$ be a local Cohen--Macaulay (CM) ring of dimension one.  It is known that $R$ has finite CM type if and only if $R$ is reduced and has bounded CM type. Here we study the one-dimensional rings of bounded but infinite CM type.  We will classify these rings up to analytic isomorphism (under the additional hypothesis that the ring contains an infinite field).  In the first section we deal with the complete case, and in the second we show that bounded CM type ascends to and descends from the completion.  In the third section we study ascent and descent in higher dimensions and prove a Brauer-Thrall theorem for excellent rings.
\end{abstract}

\maketitle

 Let $R$ be a Cohen--Macaulay (CM for short) local ring.  We say that $R$ has finite, respectively bounded, CM type provided there are only finitely many indecomposable maximal Cohen--Macaulay (MCM) $R$-modules up to isomorphism, respectively, there is a bound on the multiplicities of the indecomposable MCM $R$-modules.  The one-dimensional CM local rings of finite CM type have been completely classified \cite{Drozd-Roiter}, \cite{Green-Reiner}, \cite{Cimen:thesis}, \cite{Cimen:paper}, \cite{CWW}.  One consequence of the characterization is that a one-dimensional CM local ring $R$ has finite CM type if and only if the completion $\hat R$ has bounded CM type and is reduced.

Here we study the one-dimensional CM local rings of bounded but infinite CM type.  In \cite{Leuschke-Wiegand:hyperbrt} we classified the complete equicharacteristic hypersurface singularities of dimension one having bounded but infinite CM type:  Up to isomorphism, the only ones are $k[[X,Y]]/(Y^2)$ and  $k[[X,Y]]/(XY^2)$.  In \S1 of this paper we consider complete equicharacteristic rings that are not hypersurfaces, and we show that only one additional isomorphism type arises (Theorem~\ref{classifydim1}).  In the second section we prove that bounded CM type ascends to and descends from the completion.  Our main results in dimension one are summarized in Theorem~\ref{summary}. In the third section of the paper we study ascent and descent in higher dimensions.

\section{One-dimensional complete local rings}

We begin by quoting two results we will need from \cite{Leuschke-Wiegand:hyperbrt}.   Recall that an $R$-module $M$ is said to {\em have (constant) rank $r$} \cite{Scheja-Storch} provided $K\otimes_RM \cong K^r$, where $K$ is the total quotient ring of $R$ (obtained by inverting all nonzerodivisors).  We denote by $\nu_R(M)$ the number of generators required for $M$ as an $R$-module.

\begin{prop}{}\cite[Theorem 3.1]{Leuschke-Wiegand:hyperbrt}\label{hyperclassify}
 Let $R$ be a complete, equicharacteristic hypersurface of dimension one 
 (i.e., $R\cong k[[X,Y]]/(f)$ for some field $k$ and some 
 non-zero non-unit $f$ in the formal power series ring $k[[X,Y]]$). 
 Then $R$ has bounded but infinite CM type if and only if either 
 $R \cong k[[X,Y]]/(Y^2)$ or $R \cong k[[X,Y]]/(XY^2)$. Further, 
 if $R$ has unbounded CM type, then $R$ has, for each positive 
 integer $r$, an indecomposable MCM module of constant rank $r$.\end{prop}

\begin{prop}\label{bigranks} Let $(R,\m)$ be a one-dimensional 
local CM ring.  Assume that either
\begin{enumerate}
\item $R$ has multiplicity at least $4$, or
\item $R$ has a birational extension $S$ such that 
$\nu_R(S) = 3$ and $\nu_R(\frac{\m S}{\m}) > 1$.
\end{enumerate}
Then $R$ has, for each positive integer $r$, an 
indecomposable MCM module of constant rank $r$.
\end{prop}

While Proposition~\ref{bigranks} is not stated 
explicitly in \cite{Leuschke-Wiegand:hyperbrt}, 
it follows immediately from Lemma 2.2 and Theorem 2.3 there.  
Also, one needs (0.5) and (2.1) of \cite{CWW} to deduce the 
``Further" statement in Proposition~\ref{hyperclassify}.  
Before stating our main result (Theorem 1.5), we state two 
lemmas that will be useful here and in the next section. 
Part (1) of the first lemma is due to Bass \cite[(7.2)]{Bass:ubiquity}.  
We include the proof here, since the context in \cite{Bass:ubiquity} 
is a bit different from ours.

\begin{lemma}\label{endoring} Let $(R,\m)$ be a Gorenstein 
local ring of dimension one with total quotient ring $K$. 
Let $E = \End_R(\m) = \{\lambda \in K\ |\ \lambda \m \subseteq \m\}$.  
Assume that $E$ is local.
\begin{enumerate}
\item If $M$ is an indecomposable MCM $R$-module and 
$M\not\cong R$, then $M$ is an indecomposable
MCM $E$-module ({\it naturally}, that is, in such a way that 
the $R$-module structure induced by the natural map $R\hookrightarrow E$ 
agrees with the original structure).
\item Every indecomposable MCM $E$-module is an indecomposable MCM $R$-module.
\end{enumerate}
In particular, $R$ has finite (respectively, bounded) CM type if and only if $E$ has finite (respectively, bounded) CM type.
\end{lemma}

\begin{proof}  (1) (Bass, \cite{Bass:ubiquity}) Since there is no surjection from $M$ to $R$, we  have $M^* = \Hom_R(M,\m)$, which is an $E$-module.  Therefore $M\cong M^{**}$ is an $E$-module as well. Clearly $_EM$ is indecomposable and (since $E$ is local) maximal Cohen--Macaulay.

(2) Since $R\hookrightarrow E$ is module-finite, any MCM $E$-module $N$ is MCM as an $R$-module. Since $R\hookrightarrow E$ is birational and $_RN$ is torsion-free, any $R$-endomorphism of $N$ is automatically $E$-linear.  It follows that $N$ is indecomposable as an $R$-module.

The final statement is immediate from (1) and (2). \end{proof}

We are indebted to Tom Marley (private communication)  for showing us the following lemma, and to the anonymous referee for suggesting an improvement to the proof. Let
$\e(R)$ denote the multiplicity of the local ring $R$.

\begin{lemma}\label{tom}  Let $(R,\m,k)$ be a one-dimensional local CM ring with $k$ infinite, and suppose $\e(R) = \nu_R(\m) = 3$. Let $N$ be the nilradical of $R$.  Then:
\begin{enumerate}
\item $N^2 = 0$.
\item $\nu_R(N) \le 2$.
\item If $\nu_R(N) =2$, then $\m$ is generated by three elements $x,y,z$ such that $\m^2=\m x$ and $N=Ry+Rz$.
\item If $\nu_R(N) = 1$, then $\m$ is generated by three elements $x,y,z$ such that $\m^2=\m x$, $N = Rz$, and $yz = z^2 = 0$.
\end{enumerate}
\end{lemma}
\begin{proof}[Proof (Marley)] We note that $R$ has minimal multiplicity \cite{Abhyankar} and hence $\m$ has  reduction number 1.  Since the residue field is infinite, then, there is an element $x\in \m$ such that $\m^2 = x\m$.  Since $R$ is CM, $x$ is a nonzerodivisor.  We recall the formula
\cite[(1.1)]{Sally:1978}
\begin{equation}\label{Sallyineq}
\nu_R(J)\le \e(R) - \e(R/J)
\end{equation}
for an ideal $J$ of height $0$ in a one-dimensional CM local ring $R$. Now the image of $x$ is a  reduction element for $R/N^2$, so the right-hand side of (1.4.1), with $J = N^2$, is $\e(R/(x)) - \e(R/(Rx + N^2)$. But $N^2 \subseteq \m^2 \subset Rx$, so this expression is $0$. We conclude from (\ref{Sallyineq}) that $N^2 = 0$. Similarly, putting $J = N$, we see than $\nu_R(N) \le 2$.

For (3) and (4), we observe that $N/\m N \subseteq \m/\m^2,$ that is, minimal generators of $N$ are also minimal generators of $\m/x\m$ and of $\m$.  Indeed, since $x$ is a nonzerodivisor mod $N$, $N \cap x\m =xN$; it follows that $N\cap \m^2 = \m N$, and the map $N/\m N \to \m/\m^2$ is injective.  This proves (3).

To complete the proof of (4), let $z$ generate $N$.  Choose any $w$ such that $\m = Rx+Rw+Rz$. Now $zw = xg$ for some $g$, and (as $x$ is a nonzerodivisor) $g = bz$ for some $b\in R$. Then $z(w-bx) = 0$, and we may take $y = w-bx$.
\end{proof}

\begin{thm}\label{classifydim1}  Let $k$ be an infinite field.  
The following is a complete list, up to $k$-isomorphism, 
of the one-dimensional,  complete, equicharacteristic, 
CM local rings with bounded but infinite CM type and with residue field $k$:
\begin{enumerate}
\item $k[[X,Y]]/(Y^2)$;
\item $T:= k[[X,Y]]/(XY^2)$;
\item $E:=\End_T(\m_T)$, where $\m_T$ is the maximal ideal of $T$.  
We have a presentation $E\cong k[[X,Y,Z]]/(XY,YZ,Z^2)$.
\end{enumerate}
Moreover, if $(R,\m,k)$ is a one-dimensional, complete, 
equicharacteristic CM local ring and $R$ 
does {\it not} have bounded CM type, then $R$ has, 
for each positive integer $r$, an indecomposable MCM module 
of constant rank $r$. \end{thm}

\begin{proof}  The rings in (1) and (2) have  bounded but 
infinite CM type by Proposition~\ref{hyperclassify}.  
To show that $E$ has bounded CM type, it suffices, by 
Lemma~\ref{endoring}, to check that $E$ is local, a fact that 
will emerge in the next paragraph, where we verify the presentation of $E$ given in (3).  We routinely use decapitalization to denote specialization of variables.

The element $x+y$ is a nonzerodivisor of $T$, and the fraction $z := \frac{y^2}{x+y}$ is easily checked to be in $E:=\End_T(\m_T)$ but not in $T$.  Now $E = \Hom_T(\m_T,T)$ since $\m_T$ does not have $T$ as a direct summand, and it follows by duality over the Gorenstein ring $T$ that $E/T \cong \Ext_R^1(T/\m_T,T)\cong k$.  Therefore $E = T[z].$  Since $z^2 = \frac{y^2(x+y)^2}{(x+y)^2} \in \m_T$, $E$ is local.  One verifies the relations $xz =0, yz = y^2 = z^2$ in $E$.  Thus the map $S:=k[[X',Y',Z']] \twoheadrightarrow E$ (sending $X'\mapsto x, Y'\mapsto y,$ and $Z'\mapsto z$) induces a surjection from $A:= S/(X'Z', {Y'}^2-Y'Z', Y'Z'-{Z'}^2)$ to $E$. 
We will build an inverse map, but first we note that the change of variables
$X' = X, Y' = Y+Z, Z' = Y$ transforms $A$ to the ring $B:= k[[X,Y,Z]]/(XY,YZ,Z^2)$.  Since $B/(z)$ is one-dimensional and $z^2=0$,   $B$ and $A$ are one-dimensional. 

To build the inverse map, we note that we have $x'y'^2=0$ in $A$, so we have, at least, a map from $T$ to $A$ taking $x$ to $x'$ and $y$ to $y'$.  Next,
we note that the defining ideal $I$ of $A$ is the ideal of 
$2 \times 2$ minors of the matrix 
$\phi = \left[\smallmatrix X' & Y' \\ Y'-Z' & 0 \\ 0 & Z' \endsmallmatrix\right].$
By (the converse of) the Hilbert--Burch theorem \cite[(1.4.16)]{BH}, 
$I$ has a free resolution 
$$\CD 0 \to S^2 @>\phi>> S^3 \to I \to 0.\endCD$$  Therefore $A$ has projective dimension $2$ over $S$ and hence depth $1$.  Therefore $A$ is Cohen-Macaulay, and, since $\m_A^2 = (x'+y')\m$, it follows that $x'+y'$ is a non-zerodivisor of $A$.  Now our map $T \to A$ extends to a map $T[\frac{1}{x+y}] \to
A[\frac{1}{x'+y'}]$, and the restriction of this map provides the desired
map from $E$ to $A$.

We now know that each of the rings on our list has bounded but infinite CM type.  To show that the list is complete and to prove the ``Moreover" statement, assume now that $(R,\m,k)$ is a one-dimensional, complete, equicharacteristic CM local ring with with $k$ infinite and having infinite CM type. Suppose, moreover, that $R$ does {\it not} have indecomposable MCM modules of arbitrarily large (constant) rank. We will show that $R$ is isomorphic to one of the rings on the list.

 If $R$ is a hypersurface, Proposition~\ref{hyperclassify} tells us that $R$ is isomorphic to either $k[[X,Y]]/(Y^2)$ or $k[[X,Y]]/(XY^2)$.  Thus we assume that $\nu_R(\m) \ge 3$. But $\e(R) \le 3$ by Proposition~\ref{bigranks}.  Therefore we may assume that $\e(R) = \nu_R(\m) = 3$.   Thus we are in the situation of Lemma~\ref{tom}. Moreover, the nilradical $N$ of $R$ is non-trivial, by \cite[(0.5), (1.2)]{CWW}.

We claim that $N$ is principal.  If not, then by (3) of Lemma~\ref{tom}, we can find elements $x,y,z$ in $R$ such that
\begin{equation}\label{relationsN2gen}
\m = Rx+Ry+Rz, \ \m^2=\m x, \ {\text{and}} \ \ N=Ry+Rz
\end{equation}
Put $S := R[\frac{y}{x^2},\frac{z}{x^2}] = R + R \frac{y}{x^2} + R\frac{z}{x^2}$, and note that $\m S = \m + R\frac{y}{x} + R\frac{z}{x}$.  It is  easy to verify (by  clearing denominators) that $\{1,\frac{y}{x^2},\frac{z}{x^2}\}$ is a minimal generating set for $S$ as an $R$-module, and that the images of $\frac{y}{x}$ and $\frac{z}{x}$ form a minimal generating set  for $\frac {\m S}{\m}$. Thus we are in case (2) of Proposition~\ref{bigranks}, and our basic assumption is violated.  This proves our claim that $N$ is principal.

Using Lemma~\ref{tom}(4), we find elements $x,y,z$ in $R$ such that
\begin{equation}\label{tumrelations}
\m = Rx+Ry+Rz, \ \m^2=\m x, \ N = Rz, \ {\text{and}} \ \  yz = z^2 = 0.
\end{equation}
 Since $y^2 \in \m x \subset Rx$, we see that $R/Rx$ is a three-dimensional $k$-algebra. Further, since $\cap_n (Rx^n) = 0$, it follows that $R$ is finitely generated (and free) as a module over the discrete valuation ring $V := k[[x]]$.

We claim that $R = V + Vy + Vz$ (and therefore $\{1,y,z\}$ is a basis for $R$ as a $V$-module).  To see this, we note that $R = V[[y,z]] = V[[y]] + Vz$, since $yz = z^2 = 0$.  Let $h \in V[[y]]$, say, $h = v_0 +v_1y + v_2y^2 + \dots$, with $v_i \in V$.  For each  $n \ge 1$ write $y^{n+1} = xh_n$, with $h_n \in \m^n$. (This is possible since $\m^{n+1} = x\m^n$.)  Then $v_2y^2 + v_3y^3 + \dots = x(v_2h_1 + v_3h_2 + \dots) \in xR$.  Therefore $h\in V+ Vy + xR$, and it follows that $R = V + Vy + Vz + xR$. Our claim now follows from Nakayama's Lemma.

In order to understand the structure of $R$ we must analyze the equation that puts $y^2$ into $x\m$. Thus we write  $y^2 = x^rq$, where $r \ge 1$ and $q\in \m-\m^2$.  Write $q = \alpha x + \beta y + \gamma z$, with $\alpha,\beta, \gamma \in V$.  Since $x$ is a non-zerodivisor and $yz = z^2 = 0$, we see immediately that $\alpha = 0$.  Thus we have
\begin{equation*}
y^2 = x^r(\beta y + \gamma z)
\end{equation*}
with $\beta,\gamma\in V$; moreover, at least one of $\beta,\gamma$ must be a unit of $V$ (since $q\notin \m^2$).

We claim that $r = 1$.  For suppose $ r\ge 2$. Put $v: = \frac{y}{x^2}$ and $w:= \frac{z}{x^2}$. The relations $vw = w^2 = 0$ and $v^2 = x^{r-2} (\beta v + \gamma w)$ show that $S:= R[v,w]$ is a module-finite birational extension of $R$.  Moreover, $\nu_R(S) = 3$ and it is easy to see that $\frac{\m S}{\m}$ is minimally generated by the images of $\frac{y}{x}$ and $\frac{z}{x}$.  The desired contradiction now follows from Proposition 1.2.

Thus we have
\begin{equation}\label{y-square}
y^2 = x(\beta y + \gamma z)
\end{equation}
with $\beta,\gamma\in V$, and at least one of $\beta,\gamma$ is a unit of $V$. We will produce a hypersurface subring $A:=V[[g]]$ of $R$ such that $R = \End_A(\m_A)$.  We will then show that $A \cong k[[X,Y]]/(XY^2)$, and the proof will be complete.

{\bf Case 1}: $\beta$ is a unit. Consider the subring $A:=V[y]\subset R$.  From (\ref{y-square}) we see that $xz\in A$, and it follows easily that $z \in \End_A(\m_A)$, so $R \subset \End_A(\m_A)$.  As before (since $A$ is Gorenstein) $R = \End_A(\m_A)$.

{\bf Case 2} $\beta$ is not a unit (whence $\gamma$ is a unit). This time we put $A:=V[y+z]\subset
R$. The equation
\begin{equation*}
xy(1-\beta\gamma^{-1}) = x(y+z) - \gamma^{-1}(y+z)^2
\end{equation*}
shows that $xy \in A$.  Therefore $xz\in A$ as well, and as before we conclude that $R=\End_A(\m_A)$.

 By Lemma~\ref{endoring}, $A$ has infinite CM type but does not have indecomposable MCM modules of arbitrarily large constant rank.  Moreover, $A$ cannot have multiplicity $2$, since it has a module-finite  birational extension of multiplicity greater than $2$.  By Proposition~\ref{hyperclassify},  $A \cong k[[X,Y]]/(XY^2)$, as desired.
\end{proof}


\section{Ascent and descent in dimension one}

In this section we show that bounded CM type passes to and from the $\m$-adic completion of an
equicharacteristic one-dimensional CM local ring $(R,\m,k)$ with $k$ infinite. Contrary to the
situation in higher dimension (see Theorem~\ref{ascentgen} below), we need not assume that $R$ is
excellent with an isolated singularity.  Indeed, in dimension one this assumption would make $\hat
R$
reduced, in which case finite and bounded CM type are equivalent \cite{CWW}.    We do, however,
insist
that $k$ be infinite, in order to use the crucial fact from \S1 that failure of bounded CM type
implies
the existence of MCM modules of unbounded {\it constant} rank and also to use the explicit equations
worked out in \cite{BGS} for the indecomposable MCM modules over $T:=k[[X,Y]]/(XY^2)$.

Given a local ring $R$ with completion $\hat R$ and a finitely generated module $M$ over $\hat R$,
we
say $M$ is {\it extended} (from $R$) provided there is a finitely generated $R$-module $A$ such that
$\hat A \cong M$.  The following proposition appears, in a narrower context, in the 2002 University
of
Nebraska Ph.D. thesis of M.~Arnavut \cite{Arnavut:thesis}. The argument is adapted from
\cite[(1.5)]{Weston:1986}.

\begin{prop}\label{coextended} Let $R$ be a one-dimensional CM local ring with completion $\hat R$,
and
let $K$ be the total quotient ring of $\hat R$.
Let $M$ and $N$ be finitely generated $\hat
R$-modules
such that $K\otimes_{\hat R}M \cong K\otimes_{\hat R}N$.
Then $M$ is extended if and only if $N$ is
extended.
\end{prop}

\begin{proof} Assume $N$ is extended,
say, $N\cong \hat A$, where $_RA$ is finitely generated.
Choose
an $R$-module homomorphism $\phi:M \to N$
such that $\phi\otimes 1_K$ is an isomorphism. We obtain
an
exact sequence
\begin{equation*}\CD
0 \to V \to M @>\phi>> N @>\pi>> W \to 0,
\endCD
\end{equation*}
in which both $V$ and $W$ are torsion modules (therefore of finite
length).  Let $L = \phi(M)=\Ker(\pi)$.  Now $\Hom_{\hat R}(N, W) =
\Hom_{\hat R}(\hat A, \hat W) = (\Hom_R(A,W))\hat{} =
\Hom_R(A,W)$, since $W$ and $\Hom_R(A,W)$ have finite length.
Therefore there is a homomorphism $\rho \in \Hom_R(A,W)$ such that
$\hat \rho = \pi$.  Letting $B = \Ker(\rho)$, we see that $\hat B
\cong L$.

Next, we consider the short exact sequence
\begin{equation}\label{quartz}
0 \to V \to M \to L \to 0,
\end{equation}
viewed as an element of $\Ext_{\hat R}^1(L,V) = \Ext_R^1(B,V)$ (again, because $V$ and
$\Ext_R^1(B,V)$
have finite length). Therefore (\ref{quartz}) is the completion of a short exact sequence
\begin{equation*}
0 \to V \to C \to B \to 0.
\end{equation*}
Then $\hat C \cong M$, as desired.
\end{proof}

A finitely generated module $M$ over a Noetherian ring $R$ is said
to be {\it generically free} provided $M_P$ is $R_P$-free for each
$P\in \Ass(R)$.  For a generically free $R$-module $M$, we let
$\rank_P(M)$ denote the rank of the free $R_P$-module $M_P$, for
$P\in \Ass(R)$.

\begin{cor}\label{rankfibers} Let $(R,\m)$ be a one-dimensional CM local ring with completion $\hat
R$,
and let $M$ be a generically free $\hat R$-module.  Then $M$ is extended from $R$ if and only if
$\rank_P(M) = \rank_Q(M)$ whenever $P$ and $Q$ are minimal primes of $\hat R$ lying over the same
prime
of $R$.  In particular, every $\hat R$-module of constant rank is extended from an $R$-module
(necessarily of the same constant rank).
\end{cor}

\begin{proof} Suppose $M \cong \hat W$, and let $P$ and $Q$ be primes of $\hat R$ lying over $p\in
\Spec(R)$.  Let $r = \rank_P(M)$. We have a flat local
homomorphism $R_p \to \hat R_P$.  It follows from faithfully flat
descent \cite[(2.5.8)]{EGA4.2} that $W_p$ is $R_p$-free of rank
$r$.  From the change of rings $R_p \to \hat R_Q$ we see that
$\rank_Q(M) = r$.  This proves the ``only if'' implication and the
parenthetical remark in the last sentence of the statement.

For the converse, let $\{p_1,\dots, p_s\}$ be the minimal primes of $R$, and let $r_i =
\rank_P(M_P)$
for $P$ in the fiber over $p_i$. Let $J_1\cap\dots\cap J_s$ be a primary decomposition of $(0)$ in
$R$,
with $\sqrt{J_i} = p_i$.  (Since $R$ is CM, $\m\notin \Ass(R)$.)  Put $W =
\bigoplus_{i=1}^s(R/J_i)^{r_i}$.  Then $K\otimes_{\hat R}\hat W \cong K\otimes_R M$.  By
Proposition~\ref{coextended}, $M$ is extended.
\end{proof}

Here is our main result of this section.

\begin{thm}\label{updowndim1} Let $(R,\m,k)$ be a one-dimensional equicharacteristic CM local ring
with
completion $\hat R$.  Assume that $k$ is infinite.  Then $R$ has bounded CM type if and only if
$\hat
R$ has bounded CM type.  If $R$ has unbounded CM type, then $R$ has, for each $r$, an indecomposable
MCM module of constant rank $r$.
\end{thm}

\begin{proof} Assume that $\hat R$ does not have bounded CM type. Fix a positive integer $r$.  By
Theorem~\ref{classifydim1} we know that $\hat R$ has an indecomposable MCM module $M$ of constant
rank
$r$.  By Corollary~\ref{rankfibers} there is a finitely generated $R$-module $N$, necessarily MCM
and
with constant rank $r$, such that $\hat N \cong M$. Obviously $N$ too must be indecomposable.

Assume from now on that $\hat R$ has bounded CM type. If $\hat R$
has {\it finite} CM type, the same holds for $R$, \cite{CWW}.
Therefore we assume that $\hat R$ has infinite CM type. Then $\hat
R$ is isomorphic to one of the rings of
Theorem~\ref{classifydim1}: $k[[X,Y]]/(Y^2), T :=
k[[X,Y]]/(XY^2),$ or $E:= \End_T(\m_T).$   If $\hat R \cong
k[[X,Y]]/(Y^2)$, then  $\e(R) = 2$, and  $R$ has bounded CM type
by \cite[(2.1)]{Leuschke-Wiegand:hyperbrt}. Suppose for the moment
that we have verified bounded CM type for any local ring $S$ whose
completion is isomorphic to $E$. If, now, $\hat R \cong T$, put $S
:= \End_R(\m)$.  Then $\hat S \cong E$, whence $S$ has bounded CM
type.  Therefore so has $R$, by Lemma~\ref{endoring}.

Therefore we may assume that $\hat R = E$.  Our plan is to examine each of the indecomposable
non-free
$E$-modules and then use Proposition~\ref{coextended} and Corollary~\ref{rankfibers} to determine
exactly
which MCM $E$-modules are extended from $R$.  From now on we use the presentation for
$E=\hat R$ given in the proof of Theorem~\ref{classifydim1}: $E \cong k[[X,Y,Z]]/(XZ, Y^2-YZ, YZ-Z^2)$.  By Lemma~\ref{endoring} the indecomposable non-free MCM $E$-modules are exactly the
indecomposable non-free MCM $T$-modules, namely, the cokernels of the following matrices over $T$
(see
\cite[(4.2)]{BGS}):
\begin{equation}\label{Tmods1}
[y];\ \ [xy];\ \  [x];\ \  [y^2]
\end{equation}
\begin{equation}\label{Tmods2}
\alpha:=\left[\begin{matrix}  y & x^k\\ 0 & -y \end{matrix}\right];\ \
\beta:=\left[\begin{matrix}  xy & x^{k+1}\\ 0 & -xy \end{matrix}\right];\ \
\gamma:=\left[\begin{matrix}  xy & x^k\\ 0 & -y \end{matrix}\right];\ \
\delta:=\left[\begin{matrix}  y & x^{k+1}\\ 0 & -xy \end{matrix}\right].
\end{equation}
Let $P:=(x)$ and $Q:=(y)$ be the two minimal prime ideals of $T$.  Note that $T_P \cong k((Y))$ and
$T_Q \cong k((X))[Y]/(Y^2)$.  With the exception of $U:= \Coker [y]$ and $V:=\Coker [xy]$, each of
the
modules in (\ref{Tmods1}) and (\ref{Tmods2}) is generically free, with $\phi,
(\rank_P\Coker(\phi),\rank_Q\Coker(\phi))$ given in the following list:
\begin{equation}\label{Tmods-ranks}
[x],(1,0);\ \   [y^2],(0,1);\ \   \alpha,(0,1);\ \
\beta,(2,1);\ \  \gamma,(1,1);\ \  \delta,(1,1).
\end{equation}
Let $M$ be a MCM $\hat R$-module, and write
\begin{equation}\label{Emod}
M \cong (\oplus_{i=1}^aA_i)\oplus (\oplus_{j=1}^bB_j)\oplus(\oplus_{k=1}^cC_k)\oplus
(\oplus_{l=1}^dD_l)\oplus U^e \oplus V^f,
\end{equation}
where the $A_i,B_j,C_k,D_l$ are indecomposable generically free modules of ranks $(1,0),$ $(0,1),$
$(1,1),$ $(2,1)$ (and, again, $U = \Coker([y])$ and $V = \Coker([xy]))$.

Suppose first that $R$ is a domain.  Then $M$ is extended if and only if $b = a+d$ and $e = f= 0$.
Now
the indecomposable MCM $R$-modules are those  whose completions have $(a,b,c,d,e,f)$ minimal and
non-trivial with respect to these relations.  (We are using implicitly the fact (see
\cite[(1.2)]{Wiegand:2001} for example) that for two finitely generated $R$-modules $N_1$ and $N_2$,
$N_1$ is isomorphic to a direct summand of $N_2$ if and only if $\hat N_1$ is isomorphic to a direct
summand of  $\hat N_2$.)  The only possibilities are $(0,0,1,0,0,0)$, $(1,1,0,0,0,0)$ and
$(0,1,0,1,0,0)$, and we conclude that  the indecomposable $R$-modules have rank $1$ or $2$.

Next, suppose that $R$ is reduced but not a domain.  Then $R$ has
exactly two  minimal prime ideals, and we see from
Corollary~\ref{rankfibers} that every generically free $\hat
R$-module is extended from $R$; however, neither $U$ nor $V$ can
be a direct summand of an extended module. In this case, the
indecomposable MCM $R$-modules are generically free, with ranks
$(1,0), (0,1), (1,1)$ and $(2,1)$ at the minimal prime ideals.

Finally, we assume that $R$ is not reduced.  We must now consider the two modules $U$ and $V$ that
are
not generically free.  We will see that $U:=\Coker[y]$ is always extended and that $V:=\Coker[xy]$
is extended if and only if $R$ has two minimal prime ideals.  Note that $U \cong Txy = Exy$ (the
nilradical of $E = \hat R$), and $V \cong Ty = Ey$.

The nilradical $N$ of $R$ is of course contained in the nilradical $Exy$ of $\hat R$.  Moreover,
since
$Exy \cong E/(0:_E xy) = E/(y,z)$ is a faithful cyclic module over $E/(y,z) \cong k[[x]]$,  every
non-zero submodule of $Exy$ is isomorphic to $Exy$.  In particular, $N\hat R \cong Exy$.  This shows
that $U$ is extended.

Next we deal with $V$.  The kernel of the map $Ey \twoheadrightarrow Exy$ (multiplication by $x$) is
$Ey^2$.  Thus we have a short exact sequence
\begin{equation}\label{2.3.4}
0\to Ey^2\to V \to U\to 0.
\end{equation}
Now $Ey^2 = Ty^2 = \Coker([x])$ is generically free of rank $(1,0)$, and since the total quotient
ring
$K$ (of both $T$ and $\hat R$) is Gorenstein we see that
\begin{equation*}
KV \cong Ky^2\oplus KU.
\end{equation*}

If, now, $R$ has two minimal primes, every generically free $\hat R$-module is
extended, by Corollary~\ref{rankfibers}.  In particular, $Ey^2$ is extended, and by
Proposition~\ref{coextended} so is $V$.   Thus every indecomposable MCM $\hat R$-module is extended,
and $R$ has bounded CM type.

If, on the other hand, $R$ has just one minimal prime ideal, then the module $M$ in (\ref{Emod}) is
extended if and only if $b = a+d+f$. The $\hat R$-modules corresponding to indecomposable MCM
$R$-modules are therefore $U$, $V\oplus W$, where $W$ is some generically free module of rank
$(0,1)$,
and the modules of constant rank $1$ and $2$ described above.
\end{proof}

We conclude this section with a summary of the main results of \S\S1 and 2.

\begin{thm}\label{summary} Let $(R,\m,k)$ be an equicharacteristic one-dimensional local CM ring
with $k$ infinite.  Then $R$ has bounded but infinite CM type if and only if the completion $\hat R$
is isomorphic to one of the following:
\begin{enumerate}
\item $k[[X,Y]]/(Y^2)$;
\item $T:= k[[X,Y]]/(XY^2)$;
\item $E:=\End_T(\m_T)$, where $\m_T$ is the maximal ideal of $T$.
\end{enumerate}
Moreover, if $R$ does
{\it not} have bounded CM type, then $R$ has, for each positive integer $r$, an indecomposable MCM
module of constant rank $r$.
\end{thm}

\section{Ascent, descent, and Brauer--Thrall in higher dimensions}

In this section we study ascent and descent of bounded CM type to and from the completion in
dimension
greater than one.  We prove that bounded CM type ascends to the completion of an excellent CM local
ring with an isolated singularity.  An easy corollary of this result is a generalization of
 the Brauer--Thrall theorem of Yoshino
and Dieterich \cite[(6.4)]{Yoshino:book}.  See Theorem~\ref{BrauerThrall}.

For descent, we have a less complete picture.  We show that bounded CM type descends from the
completion of a Henselian local ring, and we investigate the case of a two-dimensional normal local
domain such that the completion is also a normal domain.  Example~\ref{nobound} indicates why
descent
is less tractable than ascent.

For Henselian rings, ascent and descent are easy:

\begin{prop}\label{updownRh} Let $R$ be a Henselian local ring with completion $\hat R$. If
$\hat R$ has  bounded CM type, then so has $R$.  Conversely, if $R$ has bounded CM type and
$\hat R$ has at most an isolated singularity, then  $\hat R$ has
bounded CM type. \end{prop}

\begin{proof} Assume $\hat R$ has bounded CM type, and let $M$ be an indecomposable MCM $R$-module.
Since $R$ is Henselian, the endomorphism ring $E:= \End_R(M)$ is local, meaning
$E/J$ is a division ring (where $J$ is the Jacobson radical of $E$).  Passing to $\hat R$, we
observe that  $\End_{\hat R}(\hat M) =
\hat R \otimes_R \End_R(M)$.  Since $\hat E/\hat J = (E/J)\hat{} = E/J$, we see that $\End_{\hat
R}(\hat M)$
is local as well, so $\hat M$ is indecomposable. Since $M$ was arbitrary, it follows that
$R$ has bounded CM type.

For the converse we use Elkik's theorem \cite[Th\`eor\'eme 3]{Elkik} on extensions
of vector bundles over Henselian pairs.  Since $\hat R$ has an isolated singularity, every MCM $\hat
R$-module $M$ is locally free on the punctured spectrum of $\hat R$, and so by Elkik's theorem is
isomorphic to $\hat N$ for some (necessarily MCM) $R$-module $N$. It follows immediately that
bounded CM type extends to $\hat R$.  \end{proof}

For ascent to the Henselization we recycle an argument from \cite{Wiegand:1998} and
\cite{Leuschke-Wiegand:2000}.  Recall \cite{Demeyer-Ingraham} that the extension $R \to \Rh$ is {\it
separable}, meaning that the sequence
\begin{equation}\label{sep}
\CD
0 @>>> J @>>> \Rh\otimes_R \Rh @>\mu>> \Rh @>>> 0,
\endCD
\end{equation}
where $\mu(u\otimes v) = uv$, is split exact as a sequence of $\Rh\otimes_R \Rh$-modules.  Tensoring
(\ref{sep}) with an arbitrary finitely generated $\Rh$-module $N$ shows that $N$ is a direct summand
of
the extended module $\Rh \otimes_R N$, where the action of $\Rh$ on $\Rh \otimes_R N$ is by change
of
rings. Write $_RN$ as a directed union of finitely generated $R$-modules $A_\alpha$.  Then, since
$N$
is a finitely generated $\Rh$-module, $N$ is a direct summand of $\Rh \otimes_R A_\alpha$ for some
$\alpha$.  Thus any finitely generated $\Rh$-module $N$ is a direct summand of an extended module.

\begin{prop}\label{recycle} Let $R$ be a CM local ring with Henselization $\Rh$.  Assume that $\Rh$
is
Gorenstein on the punctured spectrum.  If $R$ has bounded CM type, then $\Rh$ does as well.
\end{prop}

\begin{proof} Let $M$ be an indecomposable MCM $\Rh$-module. Put $d = \dim(R)$.
 Since $\Rh$ is Gorenstein on the
punctured spectrum, $M$ is a $d^{\text th}$ syzygy of some finitely generated $\Rh$-module $N$,
\cite[3.8]{EG}.  By the argument above, $N$ is a direct summand of $\Rh \otimes_R B$ for some
finitely
generated $R$-module $B$.  Letting $A$ be a $d^{\text th}$ syzygy of $_RB$, we see (using the
Krull--Schmidt Theorem over $\Rh$) that $M$ is a direct summand of $\Rh\otimes_RA$. Since $_RA$ is
MCM, we can write $A$ as a direct sum of modules $C_i$ of low multiplicity.  Using Krull--Schmidt
again, we deduce that $M$ is a direct summand of some $\Rh\otimes_RC_i$, thereby getting a bound on
the multiplicity of $M$.
 \end{proof}

\begin{thm}\label{ascentgen}  Let $R$ be an excellent CM local ring with at most an isolated
singularity.  If $R$ has bounded CM type, then the completion $\hat R$ also has bounded CM
type.\end{thm}

\begin{proof} As $R$ has geometrically regular formal fibres, $\hat R$ and the Henselization $\Rh$
both
also have isolated singularities \cite[(6.5.3)]{EGA4.2}.  By Proposition~\ref{updownRh} and
Corollary~\ref{recycle}, then, bounded CM type ascends to $\Rh$ and thence to $\hat R$.\end{proof}

Theorem~\ref{ascentgen} allows us to verify a version of the Brauer--Thrall conjecture.  The
complete
case of this theorem is due to Yoshino and Dieterich \cite[(6.4)]{Yoshino:book}.

\begin{thm}\label{BrauerThrall} Let $(R, \m, k)$ be
an excellent equicharacteristic CM local ring
with perfect residue field $k$.  Then $R$ has finite CM type if
and only if $R$ has bounded CM type and $R$ has at most an
isolated singularity.  \end{thm}

\begin{proof} If $R$ has finite CM type, then $R$ has at most an isolated singularity by
\cite{Huneke-Leuschke}, and of course $R$ has bounded CM type.  Suppose now that $R$ has bounded CM
type
and at most an isolated singularity.
According to Theorem~\ref{ascentgen}, $\hat R$ also has bounded CM type.  By the
Brauer--Thrall theorem of Yoshino and Dieterich, $\hat R$ has finite CM type.  This descends to $R$
by
\cite[(1.4)]{Wiegand:1998}, and we are done.\end{proof}

One cannot remove the hypothesis of excellence.  For example, let $S$ be any one-dimensional
analytically ramified local domain.  It is known \cite[pp. 138--139]{Matlis} that there is a
one-dimensional local domain $R$ between $S$ and its quotient field such that $\e(R) = 2$ and $\hat
R$ is not reduced.  Then $R$ has bounded but infinite CM type by \cite[(2.1),
(0.1)]{Leuschke-Wiegand:hyperbrt}, and of course $R$ has an isolated singularity.

\medskip

Proving descent of bounded CM type in general seems quite difficult.  Part of the difficulty lies in
the fact that, in general, there is no bound on the number of indecomposable MCM $\hat R$-modules
required to decompose the completion of an indecomposable MCM $R$-module.  Here is an example to
illustrate.  Recall \cite[Prop. 3]{RWW} that when $R$ and $\hat R$ are two-dimensional normal
domains, a torsion-free $\hat R$-module $M$ is extended from $R$ if and only if $[M]$ is in the
image of the
natural map on divisor class groups $\Cl(R) \to \Cl(\hat R)$.

\begin{eg}\label{nobound} Let $A$ be a complete local two-dimensional normal domain containing a
field, and assume that the divisor class group $\Cl(A)$ has an
element $\alpha$ of infinite order. (See, for example,
\cite[(3.4)]{Wiegand:2001}.) By Heitmann's theorem
\cite{Heitmann}, there is a unique factorization domain $R$
contained in $A$ such that $\hat R = A$. Choose, for each integer
$n$, a divisorial ideal $I_n$ corresponding to $n\alpha \in
\Cl(\hat R)$.  For each $n \ge 1$, let $M_n:= I_n \oplus N_n$,
where $N_n$ is the direct sum of $n$ copies of $I_{-1}$.  Then
$M_n$ has trivial divisor class and therefore is extended from $R$
by \cite[Prop. 3]{RWW}. However, no non-trivial proper direct
summand of $M_n$ has trivial divisor class, and it follows that
$M_n$ (a direct sum of $n+1$ indecomposable $\hat R$-modules) is
extended from an indecomposable MCM $R$-module.
\end{eg}

It is important to note that the example above does not give a
counterexample to descent of bounded CM type, but merely points
out one difficulty in studying descent.

\medskip

We finish with a positive result.  Recall
\cite{Geroldinger-Schneider} that the {\it Davenport constant}
$\D(G)$ of a finite abelian group $G$ is the least positive
integer $d$ such that every sequence of $d$ (not necessarily
distinct) elements of $G$ has a non-empty subsequence whose sum is
$0$. It is easy to see that $\D(G) \le |G|$, with equality if $G$
is cyclic.

\begin{prop} Let $R$ be a local normal domain of dimension two
such that $\hat R$ is also a normal domain.  Assume that the
cokernel $G$ of the natural map $\Cl(R) \to \Cl(\hat R)$ is
finite. If the ranks of the indecomposable MCM $\hat R$-modules
are bounded by $n$, then the ranks of the indecomposable MCM
$R$-modules are bounded by $m:=n\D(G)$.  In particular, if $\hat
R$ has bounded CM type, so has $R$.\end{prop}

\begin{proof} Let $M$ be an indecomposable MCM $R$-module, and put
$N:= \hat M$.  Write
\begin{equation*}
 N \cong  Z_1^{a_1} \oplus \dots \oplus Z_k^{ a_k},
\end{equation*}
where the $Z_i$ are indecomposable MCM modules, each of rank at
most $n$. It will suffice to show that $a_1+\dots+a_k \le \D(G)$.

Let $\pi:\Cl(\hat R) \to G$ be the natural map, let $\zeta_i$ be
the divisor class of $Z_i$, and let $\gamma_i = \pi(\zeta_i) \in
G$. If $a_1+\dots+a_k > \D(G)$, there is a sequence
$(b_1,\dots,b_k)\in \N^k$ such that
\begin{enumerate}
    \item $b_i \le a_i$ for all $i$;
    \item $b_j < a_j$ for some $j$;
    \item $b_k > 0$ for some $k$;
    \item $b_1\gamma_1+ \dots + b_k\gamma_k = 0$.
\end{enumerate}

Put $N' := Z_1^{b_1} \oplus \dots \oplus Z_k^{ b_k}$.  By (4), the
divisor class of the module $N'$ is in the image of the map
$\Cl(R) \to \Cl(\hat R)$.  Therefore $N'$ is extended from $R$, by
\cite[Prop. 3]{RWW}, say $N' \cong \hat M'$. Also, (1), (2) and
(3) say that $N'$ is a non-trivial proper direct summand of $N$.
By \cite[(1.2)]{Wiegand:2001}, $M'$ is a non-trivial proper direct
summand of $M$, contradiction.
\end{proof}

\providecommand{\bysame}{\leavevmode\hbox to4em{\hrulefill}\thinspace}

\end{document}